\documentclass{ifacconf}
\counterwithin*{section}{part}
\usepackage{amsmath,amsfonts,amssymb}
\usepackage{algorithmic}
\usepackage{algorithm}
\usepackage{array}
\usepackage{textcomp}
\usepackage{stfloats}
\usepackage{url}
\usepackage{verbatim}
\usepackage{multirow}
\usepackage{graphicx}
\usepackage{multirow}
\usepackage{lipsum}

\usepackage{natbib}        
\usepackage[capitalise]{cleveref}
\begin{document}
\begin{frontmatter}

\title{A PDE-Based Model of Runner Flow in Mass Participation Events} 
\thanks{This material is based upon work supported by NSF CAREER Award 2042354. The authors gratefully acknowledge Darris Blackford, Race Director, and the Columbus Nationwide Children’s Hospital Marathon for providing the data used in this study.}
\author[Second]{Alice Molinari} 
\author[First]{Brian Block} 
\author[First]{Stephanie Stockar}
\author[Second]{Simona Sacone}

\address[First]{Department of Mechanical and Aerospace Engineering and the Center for Automotive Research, The Ohio State University, 930 Kinnear Road, Columbus, OH 43212 USA.}
\address[Second]{Department of Informatics, Bioengineering, Robotics, and Systems Engineering, University of Genova, Italy.}

\begin{abstract}                
This paper presents a macroscopic approach for modeling pedestrian movement. A convection-diffusion partial differential equation (PDE) is described and used to model the flow of runners along a marathon course. A marathon is a good modeling case study for human movement because the location of runners is only known in sparse locations, so an accurate model is needed to know where people are at all times. As well, marathons offer many safety and logistical problems that can be solved by modeling human movement. Data is obtained from the 2023 Columbus Marathon, and model improvements are made to include a temporally varying convection term. Results show an increase in model accuracy which can allow for better resource allocation and safety measures.
\end{abstract}

\begin{keyword}
Modeling of pedestrian movement, logistic planning, model validation.
\end{keyword}

\end{frontmatter}

\section{Introduction}

Large pedestrian events such as concerts, tourist events, and road running races can cause severe pedestrian congestion and lead to safety problems. The Boston Marathon bomb attack in 2013 provides a striking example of how crucial the management of large-scale events is and how public safety can be significantly improved by understanding and controlling the dynamics of crowd movement \citep{CHICAGOMODEL}. As well, logistical challenges can arise in less severe scenarios such as tourist season in Cinque Terre, where large numbers of people disembark from trains at very small stations, often leading to congestion and delays \citep{[25]}.

In order to model pedestrian movement and flows in large events two types of models are used, microscopic and macroscopic. In microscopic models, the focus is on individual movement among a crowd, with the movement being modeled by an ordinary differential equations (ODEs) \citep{[2]}. One common type of microscopic pedestrian model is the cellular automaton model in which movement decisions for individuals are based on the state of neighboring cells \citep{[3]} and the speed at which pedestrians move can be modeled via a probability distribution \citep{[7]}. Another microscopic model is the social force model \citep{[8],[11]} which uses attractive and repulsive effects to determine pedestrian behavior. Social force models have also been extended to incorporate the influence of the activity pedestrians are doing, whether that is running a race, leaving a large event, or visiting a popular tourist destination \citep{[9], [10]}. While microscopic models accurately replicate pedestrian movement, they often suffer from high computational costs \citep{[2]}.

Instead, pedestrian movement can be modeled as a flow of people along a route, or numerous routes, as is done in macroscopic partial differential equation (PDE) models \citep{[1]}. Macroscopic models were first introduced in \cite{[12]} where pedestrian movement equations were derived from fluid flow. This model allowed for pedestrians to choose their own speed, but suffered from the assumption that the crowd was considered homogeneous. Another approach in \cite{[14]}, based on fluid principles, allowed pedestrians to chose their own desired speed. These fluid dynamic models, though, suffered in that they were only applicable in high density scenarios \citep{[15]}. Other macroscopic models, such as that presented in \cite{[16]}, aimed at modeling specific panic scenarios, which again were high density events.

This paper presents a novel macroscopic approach for modeling the flow of pedestrians, specifically along a marathon course. Every year in the USA, around 1000 marathons take place, some of them involving a large number of participants, such as the New York City Marathon that had 55,529 participants in 2024 \citep{nyrr2024results}, which makes marathons a great case study for model development. The model developed here is an improvement upon the model in \citep{CHICAGOMODEL}, which is a convection-diffusion PDE which captures both the forward movement of the runners as well as the dispersion of runners throughout the race. The improved model uses both time and space dependent mean speed to describe runner movement as well as develops a method for determining the dispersion of runners. Often times, data is not available for the entire marathon course, just at certain checkpoints along the course, so being able to accurately model runner flow in parts of the course that don't have data is paramount for safety and resource allocation. For example, it could be possible to predict critical times at each aid station, enabling the optimal assignment of volunteer tasks. Or, knowing when the last runners will pass through specific roads could allow for them to be reopened earlier and would provide a more accurate timetable. Additionally, choke points could be identified along the course where runners tend to crowd up, and a course reroute could be advised to alleviate this.

This paper is organized as follows. \cref{sec:model} presents the convection-diffusion model used. \cref{sec:data} showcases the available data from the 2023 Columbus Marathon used to validate the model. Then, \cref{sec:spacemeanspeed} presents the results of the base model, with a temporally constant mean speed for the convection term. Finally, \cref{sec:timemeanspeed} presents an updated calculation for the convection term that includes a variably changing mean speed to incorporate runner slow down and different classes of runners, such as elite marathoners and average, everyday runners.

\section{Convection Diffusion Pedestrian Model}\label{sec:model}
Given that a marathon course is inherently unidirectional, a one-dimensional model is used. The model is formulated as a conservation equation:
\begin{equation}
    \frac{\partial\rho(x,t)}{\partial t} = -\frac{\partial J(x,t)}{\partial x} \label{massconserv}
\end{equation}
where $x \in [0,d]$ and $t \in [0,T]$ denote distance and time, respectively, $\rho(x,t)$ is the density, and $J(x,t)$ is the flow. This model describes pedestrian movement as the combination of two different effects, convective transport and diffusion transport. Marathon dynamics can be seen as particles moving forward, however, rather than modeling the dynamics of individual particles, this model describes the dynamics of the whole runner pack. In this way, similar to the macroscale description of particle dynamics with Fokker–Planck equations \citep{CHICAGOMODEL}, the spatial and temporal changes of runner density will be considered instead.

Convection is the mass transfer as a result of the bulk motion of a fluid, or the mass transport due to the average velocity of all molecules. In the case of the marathon, runners move at a certain velocity to reach the finish line causing the convective effects. Diffusion represents the spreading of a non-uniform concentration toward a more uniform distribution over time. In a marathon this is caused by runners having different fitness levels causing them to slow down as the race progresses. The flow relationship, taking into account convective and diffusion terms, can then be written as
\begin{equation}
    J(x,t) = -D(x) \frac{\partial\rho(x,t)}{\partial x} + v(x)\rho(x,t) \label{eq:flowrhorelation}
\end{equation}
where $D(x)$ is the diffusion coefficient and $v(x)$ is the mean speed. These parameters are initially assumed to be functions of space only. 
The first term in \cref{eq:flowrhorelation}, $-D(x)\frac{\partial \rho(x,t)}{\partial x}$, represents diffusion, capturing the spreading of runners from high- to low-density regions due to variability in individual speeds. Specifically, if $\frac{\partial \rho(x,t)}{\partial x} > 0$, the diffusion term is negative, reducing the flow as density increases along the course. This reflects the natural behavior of participants that avoid moving toward crowded areas and tend to spread out where the density is lower. The extreme cases are: $D(x) =+ \infty$ where people are completely dispersed and $D(x) = 0$ where people all move in one pack at the same velocity as if they were one single body.

The second term, $v(x)\rho(x,t)$, describes the movement of people toward the finish line with a certain speed. It describes the organized motion of the group, as opposed to the random spreading captured by the diffusion term. Regions where either the speed or the density is higher contribute to a greater flow of participants. 

The complete set of equations is obtained by substituting \cref{eq:flowrhorelation} into \cref{massconserv}
\begin{equation}
    \frac{\partial \rho(x,t)}{\partial t} = \frac{\partial}{\partial x} \left( D(x)\frac{\partial \rho(x,t)}{\partial x} \right) - \frac{\partial}{\partial x} \left( v(x)\rho(x,t) \right) \label{PDE}
\end{equation}
The model is then simplified assuming that the parameters $D(x)$ and $v(x)$ are piecewise constant since data is only taken at every checkpoint along the course
\begin{equation}
    \frac{\partial \rho(x,t)}{\partial t} = D(x)\frac{\partial^2\rho(x,t)}{\partial x^2}-v(x)\frac{\partial \rho(x,t)}{\partial x}
    \label{simplPDE}
\end{equation}

For completeness, both initial and boundary conditions for \cref{PDE} must be defined. The initial condition at $t=0$, which is the beginning of the race, is defined as an empty course
\begin{equation}
    \rho (x,0) = 0
\end{equation}
The boundary condition at the beginning of the course $x=0$ is
\begin{equation}
    J_0(t) =  -D(0) \frac{\partial\rho(0,t)}{\partial x} + v(0)\rho(0,t)
    \label{feedratecondition}
\end{equation}
where $J_0(t)\in\mathbb{R}_{\geq0}$ is the feed rate. The feed rate describes how runners enter the course after the starting gun and is examined in the next section. 
At the end of the course, an open boundary condition is imposed
\begin{equation}
    \left. \frac{\partial \rho(x,t)}{\partial x} \right|_{x=L} = 0
    \label{conditiontoinf}
\end{equation}
allowing runners to freely exit once they reach the finish line.

To solve \cref{PDE}, the model is first non-dimensionalized and then discretized using the Crank–Nicolson method, a second-order finite difference scheme in time. 
The variables are non-dimensionalized to ensure appropriate scaling over the full course length. Specifically, time $t$ is scaled by the final time $T$, space $x$ by the course length $L$, and $\rho$ and $J$ by $\frac{\#\ \mathrm{runners}}{L}$ and $\frac{\#\ \mathrm{runners}}{T\cdot 60}$, respectively. 
These transformations define the dimensionless variables $\tau$, $\chi$, $\varrho$, and $\zeta$.
Substituting these new variables into \cref{simplPDE} and isolating the term of the time derivative, the non-dimensional convection-diffusion model is
\begin{equation}
    \frac{\partial \varrho}{\partial \tau}(\chi,\tau) = \delta (\chi) \frac{\partial^2\varrho}{\partial {\chi}^2}(\chi,\tau) - \omega(\chi) \frac{\partial \varrho}{\partial \chi}(\chi,\tau)
    \label{eq:nondimPDE}
\end{equation}
\begin{equation}
    -\delta(0)\frac{\partial\varrho(0,\tau)}{\partial\chi} + \omega(0)\frac{\partial\varrho}{\partial \chi}(0,\tau) = \zeta(\tau) 
    \label{eq:nondiminitialcondition}
\end{equation}
\begin{equation}
    \lim_{\chi\to 1}\frac{\partial\varrho(\chi,\tau)}{\partial\chi}=0
    \label{eq:nondimboundcondition1}
\end{equation}
\begin{equation}
    \varrho(\chi,0) = 0 
    \label{eq:nondimboundcondition2}
\end{equation}
where $\delta(\chi) = D(\chi)\frac{T}{L^2}$ is the nondimensionalized diffusion term and $\omega(\chi)=v(\chi)\frac{T}{L}$ is the nondimensionalized speed. 

To apply the Crank–Nicolson method, the continuous model must be discretized in both time and space. This involves introducing a temporal grid with time steps of size $\Delta\tau$ and a spatial grid with steps of size $\Delta\chi$.
The Crank-Nicolson method approximates the derivative of a given variable $\varrho$ with respect to $\tau$ with differences as follows: 
\begin{equation}
    \frac{\partial \varrho}{\partial\tau}(\chi,\tau+\frac{1}{2}\Delta\tau) \approx \frac{\varrho(\chi,\tau+\Delta\tau)-\varrho(\chi,\tau)}{\Delta\tau}
    \label{timederivative}
\end{equation}
This approximation is equivalent to assuming that the value of $\varrho$ at $\tau + \frac{1}{2}\Delta\tau$ corresponds to the average of its values at $\tau$ and $\tau + \Delta\tau$. In this method the interior nodes and boundary nodes are calculated differently. For interior nodes, the derivative with respect to $\chi$ is approximated via a central difference formula 
\begin{equation}
    \frac{\partial \varrho}{\partial\chi} = \frac{\varrho \left(\chi+\Delta\chi, \tau+\frac{1}{2}\Delta\tau\right)-\varrho \left(\chi-\Delta\chi, \tau+\frac{1}{2}\Delta\tau\right)}{2 \Delta\chi}
\end{equation}
Then, the second spatial derivative of $\varrho$ with respect to $\chi$ in the interior nodes is computed in $\tau+\frac{1}{2}\Delta\tau$ as
\begin{equation}
    \begin{aligned}
    \frac{\partial ^2\varrho}{\partial \chi^2} &= \frac{\varrho \left(\chi+\Delta\chi, \tau+\frac{1}{2}\Delta\tau\right)}{\Delta\chi^2} \\
    &+ \frac{\varrho \left(\chi-\Delta\chi, \tau+\frac{1}{2}\Delta\tau\right)}{\Delta\chi^2}-\frac{2\varrho \left(\chi,\tau+\frac{1}{2} \Delta\tau\right)}{\Delta\chi^2}
    \end{aligned}
\end{equation}

At the boundaries it is not possible to use the central difference formula, so a Taylor series expansion at each boundary is used. At $\chi=0$, a forward Taylor series expansion is used, which results in
\begin{equation}
    \frac{\partial \varrho}{\partial \chi}(0,\tau) = \frac{1}{\delta(0)}(\omega(0)\varrho(0,\tau)-\zeta(t))
    \label{eq:rearrangedbound}
\end{equation}
\begin{equation}
    \begin{aligned}
    \frac{\partial ^2\varrho}{\partial \chi^2}(0,\tau) &=-\frac{2}{\Delta \chi^2}\left(1+\frac{\omega(0)\Delta \chi}{\delta(0)}\right)\varrho(0,\tau)\\
    &+\frac{2}{\Delta \chi^2}\varrho(1,\tau) + \frac{2}{\Delta \chi^2}\zeta(\tau)
    \end{aligned}
\end{equation}
Combining all equations for the boundary at $\chi=0$ and computing some substitutions, the Crank-Nicholson update is
\begin{equation}
    \begin{aligned}
    &\varrho(0,\tau+\Delta \tau)\left[\frac{1}{\Delta \tau}+\left(1+\frac{\omega(0)\Delta \chi}
    {\delta(0)}\right)\frac{\delta(0)}{\Delta \chi^2}+\frac{\omega(0)^2}{2\delta(0)}\right]\\
    &+\varrho(\Delta\chi,\tau+\Delta \tau)\left[-\frac{\delta(\Delta\chi)}{\Delta \chi^2}\right] = \\
    &\varrho(0,\tau)\left[\frac{1}{\Delta \tau}+\left(1+\frac{\omega(0)\Delta \chi}{\delta(0)}\right)\frac{\delta(0)}{\Delta \chi^2}+\frac{\omega(0)^2}{2 \delta(0)}\right]\\
    &+\varrho(\Delta\chi,t)\left[-\frac{\delta(\Delta\chi)}{\Delta \chi^2}\right]+\zeta(t)\left[\frac{2 \delta(0)}{\Delta \chi^2}+\frac{\omega(0)}{\delta(0)}\right]
    \end{aligned}
\end{equation}
where all terms for $\tau+\Delta\tau$ are on the left-hand side. Following a similar approach, it is now possible to obtain the equation for $\varrho(1,\tau)$, the boundary condition at the end of the course, but this time using a backward Taylor series expansion. The Crank-Nicholson update at $\chi=1$ is thus
\begin{equation}
    \begin{aligned}
        &\left(\frac{1}{\Delta \tau} +\frac{\delta(1)}{\Delta \chi ^2}\right) \varrho(1, \tau + \Delta \tau) + \\
        &\left(-\frac{\delta(1)}{\Delta \chi ^2}\right)\varrho(1-\Delta \chi, \tau + \Delta \tau)\\
        &= \left(\frac{1}{\Delta \tau} -\frac{\delta(1)}{\Delta \chi ^2}\right) \varrho(1, \tau) + \left(-\frac{\delta(1)}{\Delta \chi ^2}\right)\varrho(1-\Delta \chi, \tau)
    \end{aligned}
\end{equation}

After discretizing the interior nodes and applying the boundary conditions, the resulting equations are assembled into a tridiagonal matrix system. This system is solved using a direct Gauss–Jordan elimination method, which applies forward elimination to obtain an upper triangular system, followed by back substitution to solve from the last equation to the first.

\section{Data Analysis}\label{sec:data}
\begin{figure}
    \centering
    \includegraphics[width=1\linewidth]{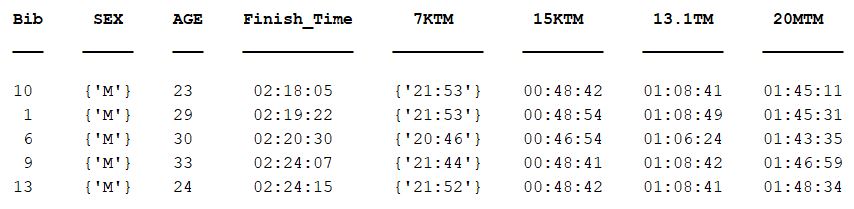}
    \caption{Sample of the Columbus Marathon datasets with runner times at multiple checkpoints.}
    \label{fig:HeadOfDataset}
\end{figure}

The data provided comes from the 2023 Columbus Marathon \citep{cbus2023results}. The data includes biographical data as well as race specific data such as chip start time, time to cross race checkpoints, and finish time for the 2967 runners. The chip start time measures when a runner crosses the start line and the race checkpoints occur at 7 km, 15 km, 13.1 m (21.1 km), and 20 m (32.2 km). An example set of data is shown in \cref{fig:HeadOfDataset}

\begin{figure}
    \centering
    \includegraphics[width=1\linewidth]{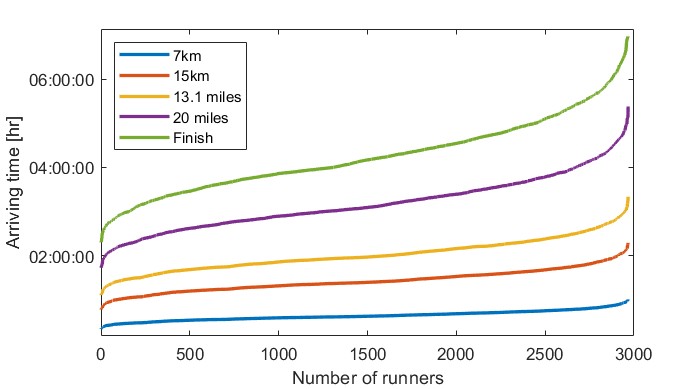}
    \caption{Runner arrival times at each checkpoint.}
    \label{fig:arrivingorder}
\end{figure}

The time at which each runner arrives at the race checkpoints is shown in \cref{fig:arrivingorder} which presents five lines, one for each checkpoint and the finish line. The fastest runners arrive at the first check-point in less than 30 minutes and finish the race in less than 2.5 hours, while it takes the slowest runners close to 7 hours to finish. At 7 km, the gap between the first runner and the last one is about forty minutes, while at the finish line this value increases to more than four hours. From the halfway point to 20 miles, the slower runners start to fatigue, a phenomenon known as 'hitting the wall' in which runners begin to exhaust physiologic carbohydrate reserves and start to slow down \cite{[24]}. This behavior is an important factor for the model to capture.

\begin{figure}
    \centering
    \includegraphics[width=1\linewidth]{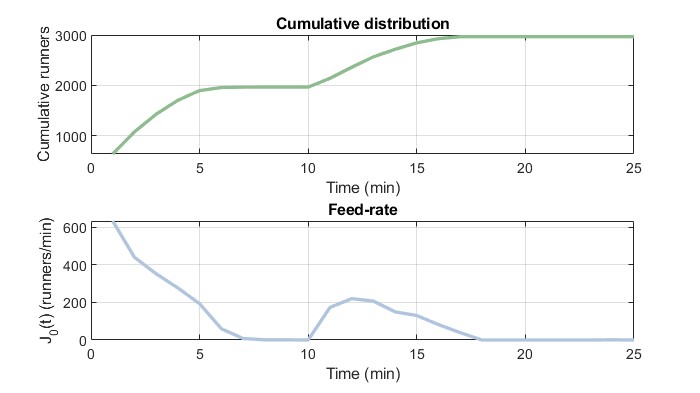}
    \caption{Cumulative number of runners and corresponding feed rate at the race start.}
    \label{fig:FeedRate}
\end{figure}

The flow and density are calculated from the provided data as
\begin{equation}
    J(x,t) = \frac{N(x,t;\Delta t)}{\Delta t},\qquad \rho(x,t) = \frac{N(x,t; \Delta x)}{\Delta x}
\end{equation}
where $J(x,t)$ is the flow of runners crossing position $x$ at time $t$, $N(x,t;\Delta t)$ is the number of runners crossing $x$ within the interval time $[t,t+\Delta t]$, and $N(x,t; \Delta x)$ represents the number of runners that move in the space interval $[x,x+\Delta x]$ at time $t$. The space interval $\Delta x$ is not constant, but is calculated as the difference between consecutive checkpoints. The boundary condition at $x=0$, or the feed rate  $J_0(t)$, is shown in \cref{fig:FeedRate}. \cref{fig:FeedRate} also shows the cumulative number of runners entering the course. The line in the top plot assumes an S-shape, creating two waves, separated by 3 minutes. This is done by the race organizers to reduce congestion along the course. The first wave corresponds to about two-thirds of the total number of runners. All runners enter the race in less than 20 minutes.

\section{Base Model Results}\label{sec:spacemeanspeed}
The mean speed of the runner pack is an important parameter of the model which represents convection transport.  The mean speed, $v(x)$ is first calculated as only a function of space. The speed of each runner at each checkpoint is calculated as
\begin{equation}
    v_{ijn} = \frac{d_{jn}-d_{in}}{t_{jn}-t_{in}}
    \label{eq:velocity}
\end{equation}
where the index $i$ represents a checkpoint and $j$ represents the following one, and the index $n$ indicates the number of the runner. In \cref{eq:velocity}, $d_{jn}$ is the distance in kilometers that corresponds to the checkpoint $j$, $d_{in}$ is the distance in kilometers that corresponds to the checkpoint $i$, $t_{jn}$ is the time after which the runner $n$ is at the checkpoint $j$, $t_{in}$ is the time after which the runner $n$ is at the checkpoint $i$, and $v_{ijn}$ is the speed of the runner $n$ between $i$ and $j$.

\begin{figure}
    \centering
    \includegraphics[width=1\linewidth]{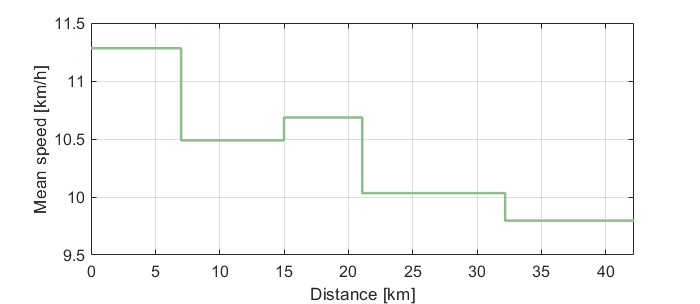}
    \caption{Mean runner speed as a function of course distance.}
    \label{fig:MeanSpeed}
\end{figure}

After calculating all the velocities of all runners, the average over all runners is computed, and a vector of five points is obtained, one for each spatial interval. \cref{fig:MeanSpeed} shows the mean speed of runners as a function of space. As expected, the mean speed is a mostly decreasing function of space, except for the interval between 15 km and 13.1 miles. This is due to marathon runners being influenced by the faster half marathon runners who are finishing, since both races occur on the same course.

\begin{figure}
    \centering
    \includegraphics[width=1\linewidth]{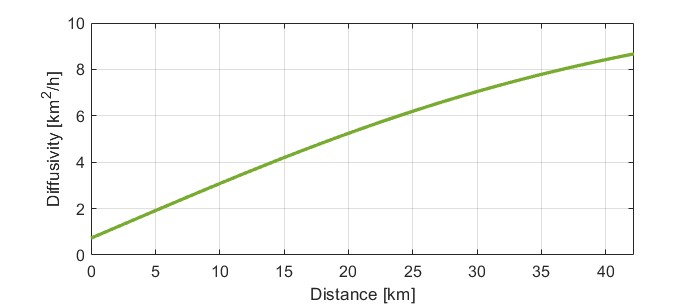}
    \caption{Diffusivity as a function of course distance.}
    \label{fig:Diffusivity}
\end{figure}

The diffusion coefficient $D(x)$ is the second parameter of the model that is used to represent diffusion transport. As with the mean speed, this parameter is a function of space, but conversely to the mean speed it is expected to be increasing in space because as the marathon goes on the runners spread more along the race course due to their different velocities and aerobic capacities. The diffusivity is described by the following equation
\begin{equation}
    D(x) = \frac{v(x)\rho(x,t)-J(x,t)}{\frac{\partial\rho(x,t)}{\partial x}}
    \label{eq:diffusivity}
\end{equation}
It is not possible to find the diffusion directly from \cref{eq:diffusivity} with the available data; therefore, a parameter identification algorithm is used which finds the diffusivity that minimizes the error between the real flow and the flow generated with the model. The algorithm used is called Particle Swarm Optimization (PSO), as described in \citep{[22]}. The function used for diffusivity is a replication of the diffusivity shown in \citep{CHICAGOMODEL}, representing the diffusion coefficient of the Chicago Marathon. The function is as follows:
\begin{equation}
    D(x) = p_1 + p_2\mathrm{tanh}(p_3(x-p_4))
    \label{replication}
\end{equation}
Here $p_1$, $p_2$, $p_3$, and $p_4$ are the parameters to be found. The PSO algorithm aims to minimize the root mean square error (RMSE)
\begin{equation}
    \mathrm{RMSE} = \sqrt{\frac{1}{T}\sum_{j=1}^{T}e(t_j)^2}
\end{equation}
where
\begin{equation}
    e_i(t_j) = J_i(t_j) -\hat{J}_i(t_j)
\end{equation}
and the index $i$ represents the different checkpoints and $j$ represents the different instants in time. $J_i(t_j)$ is the real flow that crosses the checkpoint $i$ at the instant $j$ and $\hat{J}_i(t_j)$ is the corresponding estimated flow from the model. This results in \cref{fig:Diffusivity} where, as expected, $D(x)$ is an increasing function of space, which is consistent with what has been observed in the previous analysis. This is due to the fact that runners disperse along the race course.

\begin{figure}
    \centering
    \includegraphics[width=1\linewidth]{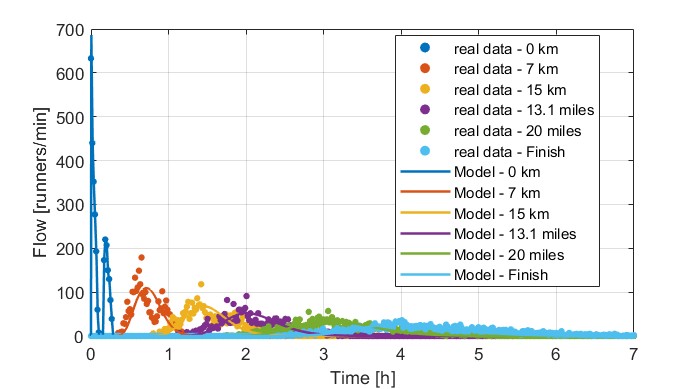}
    \caption{Measured and modeled flow rates at each checkpoint using the constant speed model.}
    \label{fig:FlowComparison}
\end{figure}

The comparison of flow rates at each checkpoint and the finish line for the real data and model is shown in \cref{fig:FlowComparison}. The model captures the main trends in the data, including the temporal spreading of flow at each checkpoint as runners disperse over time. Specifically, the curves tend to spread over time at each checkpoint due to the fact that at the beginning the participants run relatively together until they begin to spread out as they tire. The peaks indicate the average arrival time of the main group of runners at each checkpoint. This time is important for race organizers as it lets them know when the majority of runners have crossed a checkpoint so that they can better allocate volunteers or begin to open up closed roads along the route. 

\cref{tab:Results(v)} compares average checkpoint times from data and model. The values in the first column are calculated using the arithmetic mean:
\begin{equation}
    \mathrm{Mean\,Time}_i = \frac{1}{N}\sum_{j=1}^{N}t_{ji} \quad \forall i 
\end{equation}
where $i$ is the index of the checkpoint, $j$ is the index of the runner, $N$ is the total number of runners, and $t_{ji}$ is the time after the gunshot at which the runner $j$ reaches the checkpoint $i$. The second column is computed from the model results using a weighted average, where the flow acts as a weight
\begin{equation}
    \mathrm{Mean\,Time\,Model}_i = \frac{\sum_{t=1}^{T}J_{it} \cdot t_t}{\sum_{t=1}^{T}J_{it}} \quad \forall i
\end{equation}
where $t$ is the index of the discrete time step, $T$ is the total number of time steps in the model, $J_{it}$ is the modeled flow at checkpoint $i$ and at time $t$, and $t_t$ is the time elapsed since the gunshot at step $t$. The percentage difference is then computed as:
\begin{equation}
    \Delta \%_i = \frac{\mathrm{Mean\,Time}_i-\mathrm{Mean\,Time\,Model}_i}{\mathrm{Mean\,Time}_i} \cdot 100
\end{equation}
The model overestimates the time to cross the first checkpoint as shown in \cref{tab:Results(v)}, but in the back half of the race the error grows to an almost 4.4\% underestimation of the finishing time. 
\begin{table}[thpb]
    \centering
    \begin{center}
    \caption{Mean times at each checkpoint.}
    \label{tab:Results(v)}
    \begin{tabular}{l|c|c|c|c}
    \hline \hline
    Checkpoint  & Real Data [hr] & Model [hr] & $\Delta \%$ & RMSE\\\hline
    7 km        & 0.7339         & 0.7503     & -2.2339     & 11.27\\
    15 km       & 1.5211         & 1.5174     & 0.2460      & 6.57\\
    13.1 miles  & 2.1113         & 2.1061     & 0.2460      & 6.11\\
    20 miles    & 3.2704         & 3.2043     & 2.0203      & 4.61\\
    Finish      & 4.3462         & 4.1562     & 4.3714      & 3.94\\
    \hline \hline
    \end{tabular}
    \end{center}
\end{table}

\section{Improved Model Results}\label{sec:timemeanspeed}
In this section, instead of a spatially varying mean speed, $v(x)$, the model is extended to incorporate also a temporally varying mean speed, $v(x,t)$. Intuitively, mean speed varies over both time and distance as runners passing through earlier in the race will be faster and shouldn't influence the mean speed of the following runners. This effect is know as \emph{hitting the wall} \citep{Buman_Brewer_Cornelius_2009,[24]} and is caused by runners' previous training and nutrition. In the former approach in the previous section, the mean speed was averaged over all runners, but this results in an underestimation of the lead runners and an overestimation of the speed of the runners at the back.

In order to implement the new mean speed function, that is $v(x,t)$, which varies over space and time, a time interval must be defined which in this case is $\Delta t= 1$ h. Then, the approach consists of ordering from the fastest runner to the slowest each column of the dataset containing the time after which each runner reaches each checkpoint. Once the columns are sorted, runners are separated into groups based on the hour in which they reach each checkpoint. The mean speed is determined by computing the average of the velocities of the runners in each group.

\begin{figure}
    \centering
    \includegraphics[width=0.75\linewidth]{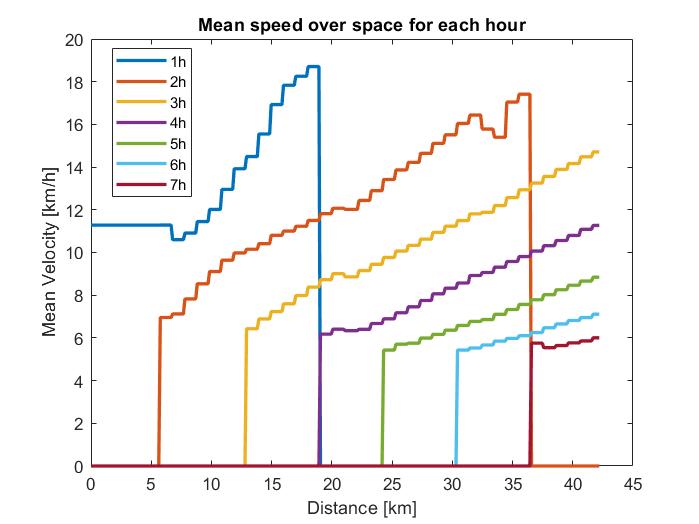}
    \caption{Mean runner speed as a function of space and time.}
    \label{fig:MeanSpeed_space}
\end{figure}

\cref{fig:MeanSpeed_space} shows the new mean speed function, where the x-axis is the distance of the race in kilometers, and the y-axis is still the mean speed in km/h. This graph shows seven different curves, one for each hour of the race. The plots go to zero where no one is running. For example, no runner runs more than about 20 km in the first hour. The second hour is the most spread function, where some people have already ran 35 kilometers, and some others are still running at the tenth kilometer. 

\begin{figure}
    \centering
    \includegraphics[width=1\linewidth]{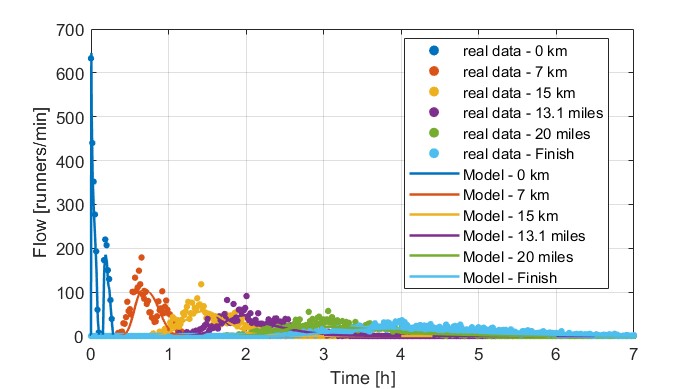}
    \caption{Measured and modeled flow rates at each checkpoint using the time-varying speed model.}
    \label{fig:Flow_comparison_x_t}
\end{figure}

Earlier runners exhibit higher mean speeds, while later groups reflect fatigue effects. In addition, fast runners finish the race within the first 2.5 hours, so they are no longer considered in the computation of the mean speed. Once the new mean speed is found, it is possible to use it as an input to the model, where at each hour a different mean speed vector is chosen. \cref{fig:Flow_comparison_x_t} shows the resulting modeled flow rate at each checkpoint, compared against the data. As in \cref{fig:FlowComparison}, the model represents the data well.

As shown in \cref{tab:Results(v_x_t)}, the improved model provides a more accurate estimate of the average runner. In the first checkpoint the error is decreased by roughly 2\% and the average finishing time error is also decreased by 2\%. For the other checkpoints, the new model estimates the average crossing time almost exactly.
\begin{table}[thpb]
    \centering
    \begin{center}
    \caption{Time comparison at every checkpoint using variable mean speed.}
    \label{tab:Results(v_x_t)}
    \begin{tabular}{l|c|c|c|c}
    \hline \hline
    Checkpoint  & Real Data [hr] & Model [hr] & $\Delta \%$ & RMSE \\\hline
    7 km        & 0.7339         & 0.7296     & 0.5889      & 12.31\\
    15 km       & 1.5211         & 1.5120     & 0.6012      & 7.01\\
    13.1 miles  & 2.1113         & 2.1104     & 0.0408      & 6.19\\
    20 miles    & 3.2704         & 3.2679     & 0.0749      & 5.05\\
    Finish      & 4.3462         & 4.2433     & 2.3670      & 4.61\\
    \hline \hline
    \end{tabular}
    \end{center}
\end{table}

\section{Conclusion}\label{sec:conc}
This paper presented a way of modeling flow during a mass event, such as a city marathon. First, a PDE model of pedestrian flow was presented and then real world data was taken from the Columbus Marathon and analyzed. This data was then used to validate the base model. Finally, model improvements were made by defining a spatially and temporally varying mean speed to better capture the characteristics of the event, such as runners slowing down as they get towards the end of the race. The improvements to the mean speed calculation were able to better model the movement of runners along the course, reducing average time to reach checkpoints along the course by up to 2\%.

Future work will look at using this model to reduce areas of high runner density along the course, to analyze the impact of changing how runners start the race, and to better plan the assignment of aid stations and resources along the course. As well, this approach will be extended to model other large events such as concerts, fairs, parades, and festivals such that organizers can create better plans for crowd control and safe evacuation.

\bibliography{bibliografia}

@article{[1],
author = {Shiwakoti, Nirajan and Sarvi, M. and Rose, G.},
year = {2008},
pages = {457-473},
title = {Modelling pedestrian behaviour under emergency conditions: State-of-the-art and future directions},
journal = {31st Australasian Transport Research Forum, ATRF 2008}
}

@article{[2],
title = {State-of-the-art crowd motion simulation models},
journal = {Transportation Research Part C: Emerging Technologies},
volume = {37},
pages = {193-209},
year = {2013},
author = {Dorine C. Duives and Winnie Daamen and Serge P. Hoogendoorn},
}

@article{[3],
author = {Victor J. Blue and Jeffrey L. Adler},
title ={Emergent Fundamental Pedestrian Flows from Cellular Automata Microsimulation},
journal = {Transportation Research Record},
volume = {1644},
number = {1},
pages = {29-36},
year = {1998},
}

@article{[7],
    author = {Song and Weiguo and Yu and Yanfei and Fan and Weicheng and Zhang and Heping} ,
    title = "A cellular automata evacuation model considering friction and repulsion",
    journal = "Science in China Ser. E Engineering \& Materials Science" ,
    year = "2005",
    volume = "48",
    pages = "403-413"
}

@article{[8],
 author = {Dirk Helbing and Lubos Buzna and Anders Johansson and Torsten Werner},
 journal = {Transportation Science},
 number = {1},
 pages = {1--24},
 publisher = {INFORMS},
 title = {Self-Organized Pedestrian Crowd Dynamics: Experiments, Simulations, and Design Solutions},
 urldate = {2025-04-16},
 volume = {39},
 year = {2005}
}

@article{[9],
title = {The Nomad Model: Theory, Developments and Applications},
journal = {Transportation Research Procedia},
volume = {2},
pages = {462-467},
year = {2014},
author = {Mario Campanella and Serge Hoogendoorn and Winnie Daamen},
}

@article{[10],
title = {Pedestrian route-choice and activity scheduling theory and models},
journal = {Transportation Research Part B: Methodological},
volume = {38},
number = {2},
pages = {169-190},
year = {2004},
author = {S.P. Hoogendoorn and P.H.L. Bovy},
}

@article{[11],
  title = {Social force model for pedestrian dynamics},
  author = {Helbing, Dirk and Moln\'ar, P\'eter},
  journal = {Phys. Rev. E},
  volume = {51},
  pages = {4282--4286},
  year = {1995},
}

@article{[12],
    author = {Henderson, L. F.},
    title = {The Statistics of Crowd Fluids},
    journal = {Nature},
    year = {1971},
    volume = {229},
    pages = {381--383},
}

@article{[14],
author = {Helbing, Dirk},
year = {1998},
pages = {391--415},
title = {A Fluid Dynamic Model for the Movement of Pedestrians},
volume = {6},
journal = {Complex Systems}
}

@article{[15],
title = {A continuum theory for the flow of pedestrians},
journal = {Transportation Research Part B: Methodological},
volume = {36},
number = {6},
pages = {507-535},
year = {2002},
author = {Roger L. Hughes},
}

@article{[16],
author = {Colombo, Rinaldo and Rosini, Massimiliano},
year = {2005},
pages = {1553-1567},
title = {Pedestrian flows and non-classical shocks},
volume = {28},
journal = {Math. Methods Appl. Sci.},
}

@article{CHICAGOMODEL,
title = {Convection–diffusion effects in marathon race dynamics},
journal = {Physica A: Statistical Mechanics and its Applications},
volume = {393},
pages = {498-507},
year = {2014},
author = {E. Rodriguez and G. Espinosa-Paredes and J. Alvarez-Ramirez},
}

@inproceedings{[22],
author = {JE, Kennedy and Eberhart},
booktitle = {Proceedings of the IEEE International Conference on Neural Networks},
year = {1995},
pages = {1942–1945},
title = {Particle Swarm Optimization},
volume = {4}
}

@misc{nyrr2024results,
  author       = {{New York Road Runners}},
  title        = {TCS New York City Marathon 2024 Results},
  year         = {2024},
  howpublished = {\url{https://results.nyrr.org/event/M2024/finishers}},
}

@misc{cbus2023results,
  author       = {{{MTEC} Results}},
  title        = {Nationwide Children's Hospital Columbus Marathon 2023 Results},
  year         = {2023},
  howpublished ={\url{www.mtecresults.com/event/show/5003/2023_Columbus_Marathon}},
}

@article{[24],
  author       = {Rapoport, Benjamin I.},
  title        = {Metabolic Factors Limiting Performance in Marathon Runners},
  journal      = {PLoS Computational Biology},
  year         = {2010},
  volume       = {6},
  number       = {10},
}

@article{[25],
author = {del Cacho Estil-les, María Asunción and Bersani, Chiara and Sacile, Roberto and Zero, Enrico},
year = {2023},
pages = {89152 - 89163},
title = {Optimal Travel Planning of Short Stays in Mass Tourist Destinations},
volume = {11},
journal = {IEEE Access},
}

@article{Buman_Brewer_Cornelius_2009, 
    title={A discrete-time hazard model of hitting the wall in recreational marathon runners}, 
    volume={10}, 
    ISSN={1469-0292}, 
    journal={Psychology of Sport and Exercise}, 
    author={Buman, Matthew P. and Brewer, Britton W. and Cornelius, Allen E.}, 
    year={2009}, 
    month=nov, 
    pages={662–666}, 
    language={en} 
}
\end{document}